\documentclass[11pt]{amsart}
\usepackage{amscd}
\usepackage{amssymb}

\newtheorem{thm}{Theorem}[section]
\newtheorem{lem}[thm]{Lemma}

\theoremstyle{definition}

\renewcommand{\thecase}{}
\newtheorem{conj}[thm]{Conjecture}

\renewcommand{\thestep}{}

\theoremstyle{remark}

\makeatletter
\def\alphenumi{
  \def\theenumi{\alph{enumi}}
  \def\p@enumi{\theenumi}
  \def\labelenumi{(\@alph\c@enumi)}}
\makeatother

\makeatletter
\def\thecase{\@arabic\c@case}
\makeatother

\numberwithin{equation}{section}

\makeatletter
\def\thestep{\@arabic\c@step}
\makeatother

\newenvironment{pf*}[1]{\begin{proof}[#1]}{\end{proof}}

\newcommand\barM{{\bar{M}}}

\newcommand\barV{{\bar{V}}}
\newcommand\barW{{\bar{W}}}

\newcommand\AAA{\mathbb{A}}

\newcommand\CC{\mathbb{C}}

\newcommand\RR{\mathbb{R}}

\newcommand\ZZ{\mathbb{Z}}

\newcommand\bD{{\mathbf{D}}}

\newcommand\bL{{\mathbf{L}}}

\newcommand\bS{{\mathbf{S}}}

\newcommand\bW{{\mathbf{W}}}

\newcommand\half{{\textstyle{\frac{1}{2}}}}

\newcommand\quarter{{\textstyle{\frac{1}{4}}}}

\newcommand\fs{{\mathfrak{s}}}

\newcommand\La{\Lambda}

\newcommand\si{\sigma}
\newcommand\Si{\Sigma}

\newcommand\su{{\mathfrak{s}\mathfrak{u}}}

\newcommand\PU{\operatorname{PU}}

\newcommand\SO{\operatorname{SO}}

\newcommand\asd{{\operatorname{asd}}}

\newcommand\Ind{\operatorname{Index}}

\newcommand\SW{SW}
\newcommand\Sym{\operatorname{Sym}}

\newcommand\even{{\mathrm{even}}}

\newcommand\spinc{\text{$\text{spin}^c$ }}
\newcommand\Spinc{\text{$\mathrm{Spin}^c$}}

\newcommand\sO{{\mathcal{O}}}

\begin{document}
\title[$\PU(2)$ monopoles and a conjecture of Mari\~no, Moore, and Peradze]
{\boldmath{$\PU(2)$} monopoles and a conjecture of Mari\~no, Moore, and Peradze
}
\author[P. M. N. Feehan]{P. M. N. Feehan}
\address{Department of Mathematics\\
Ohio State University\\
Columbus, OH 43210}
\curraddr{School of Mathematics\\
Institute for Advanced Study\\
Olden Lane\\
Princeton, NJ 08540}
\email{feehan@math.ohio-state.edu and feehan@math.ias.edu}
\author[P. B. Kronheimer]{P. B. Kronheimer}
\address{Department of Mathematics\\
Harvard University\\
1 Oxford Street\\
Cambridge, MA 02138}
\email{kronheim@math.harvard.edu}
\author[T. G. Leness]{T. G. Leness}
\address{Department of Mathematics\\
Florida International University\\
Miami, FL 33199}
\email{lenesst@fiu.edu}
\author[T. S. Mrowka]{T. S. Mrowka}
\address{Department of Mathematics\\
rm 2-367\\
Massachusetts Institute of Technology\\
Cambridge, MA 02139}
\email{mrowka@math.mit.edu}
\thanks{PMNF was supported by NSF grant number
DMS-9704174 and, through the Institute for Advanced Study, by NSF grant
number DMS-9729992; PBK was supported by NSF grant number DMS-9531964;
TSM was supported by NSF grant number DMS-9796248}



\maketitle


\section{Introduction}
\label{sec:Intro}
The purpose of this note is to show that some of the recent results of
Mari\~no, Moore, and Peradze \cite{MMPhep}, \cite{MMPdg}
can be understood in a simple and direct way
via a mechanism pointed out in \cite{FLGeorgia},
\cite{FL2}, \cite{FL2b}, using the $\PU(2)$-monopole cobordism of
Pidstrigach and Tyurin \cite{PTLocal}.

Throughout this paper, let $X$ denote an oriented smooth four-manifold with
$b^{+}_{2}(X)>1$ and with $b^{1}(X)=0$.
Once an orientation  of
$H^{2}_{+}(X)$ is chosen, we can define
the Seiberg-Witten
invariants which we view as a function
$$
SW_{X}:\Spinc(X) \to \ZZ,
$$
where $\Spinc(X)$ denotes the set of isomorphism classes of \spinc
structures on $X$. Let $S=S(X)\subset \Spinc(X)$ be the support of
$SW_{X}$.  A cohomology class $K\in H^2(X;\ZZ)$ is called an {\it SW-basic
class} if $K =c_{1}(\fs)$ for some $\fs \in S$. Let $B=B(X)\subset
H^2(X;\ZZ)$ be the set of SW-basic classes.  A four-manifold $X$ is
said to be of SW-simple type if for all $\fs\in S$
\begin{equation}
\label{eq:SimpleType}
c_{1}(\fs)^{2}=2\chi(X)+3\sigma(X).
\end{equation}
Given an integral two-dimensional cohomology class $w$, we combine the
Seiberg-Witten invariants into an analytic
function of $h\in H_{2}(X;\RR)$ by the formula
$$
\bS\bW_X^{w}(h) =\sum_{\fs \in
S(X)}(-1)^{\half(w^{2}+c_{1}(\fs)\cdot w)}
SW_X(\fs)e^{\langle c_{1}(\fs),h\rangle}.
$$
Let $B^{\bot} \subset H^{2}(X,\ZZ)$ be the orthogonal complement
of the basic classes. We say that a four-manifold is {\it abundant} if
the intersection form on $B^{\bot}$ contains a hyperbolic sublattice.
Define the characteristic number
\begin{equation}
\label{eq:cX}
c(X)=-\frac{1}{4}(7\chi(X)+11\sigma(X)).
\end{equation}
If $X$ is a complex surface then $c(X)=\chi_h(X)-c_{1}^{2}(X)$, where
$\chi_h(X)=\chi(\sO_X)$. Our main theorem is
\begin{thm}\label{main}
	Suppose that $X$ is abundant and of SW-simple type and that
	Conjecture~\ref{conj:Divisibility} holds for $X$. Then 
	either $c(X)-3< 0$, or for any integral lift $w$ of $w_{2}(X)$,
	the series $\bS\bW^{w}_{X}$ vanishes to order $c(X)-2$ at $h=0$.
\end{thm}

As explained in \cite[\S 8.1]{MMPhep}, this result constrains the homotopy
type\footnote{A similar constraint was conjectured earlier by Fintushel and
Stern (private communication; see also \cite{FSOneBasicClass}).}
of $X$ in terms of the number $b(X)$ of elements in $B/\{\pm 1\}$:

\begin{thm}
\label{GeneralizedNoether}
\cite[Theorem 8.1.1]{MMPhep}
Let $X$ be a closed, oriented, smooth four-manifold with $b_2^+(X)>1$. If
$b(X)>0$ and $X$ obeys the conclusion of Theorem \ref{main}, then 
\begin{equation}
\label{eq:NumberBasicClasses}
b(X) > c(X)/2.
\end{equation}
\end{thm}

The examples of Fintushel and Stern \cite{FSOneBasicClass} imply that the
bound of Theorem~\ref{main} is sharp: for every point along the line
$c(X) = \text{constant}\ge 2 $ in the $(c_1^2,\chi_h)$ plane
there exists a four-manifold $X$ with $\bS\bW^{w}_{X}$
vanishing to order $c(X)-2$.  Here
$c_1^2(X) = 2\chi(X)+3\sigma(X)$, and
$\chi_h(X) = \quarter(\chi(X)+\sigma(X))$. It is interesting to note that
the slope of the line in the $(c_1^2,\chi_h)$ plane implied by inequality
\eqref{eq:NumberBasicClasses}, namely
\begin{equation}
\label{eq:GeneralizedNoether}
c_1^2(X) \geq \chi_h(X) - 2b(X) - 1,
\end{equation}
does not coincide with the classical Noether or Bogomolov-Miyaoka-Yau lines for
complex surfaces.

The vanishing condition in the conclusion of Theorem~\ref{main} is the
statement that $X$ has ``superconformal simple type'' in the terminology of
\cite{MMPhep}, where it is further conjectured (Conjecture 7.8.1) that all
four-manifolds of SW-simple type have this property.  Theorem~\ref{main}
therefore reduces the conjecture of \cite{MMPhep} for abundant manifolds to
the technical Conjecture \ref{conj:Divisibility}.  The latter conjecture is
the assertion that ideal Seiberg-Witten moduli spaces
$M_\fs^{sw}\times\Sym^\ell(X)$ in any stratum of the compactified moduli
space of $\PU(2)$ monopoles $\barM_{W,E}$ make no contribution to the
Donaldson invariants defined by $\barM_E^{\asd}\subset \barM_{W,E}$ if
their associated Seiberg-Witten invariants vanish. This in turn is a simple
consequence of a conjecture, attributed to Pigstrigach and Tyurin
\cite[Conjecture 4.1]{FLGeorgia}, that pairings of Donaldson-type
cohomology classes with links of ideal Seiberg-Witten moduli spaces in
$\barM_{W,E}$ are given by the product of $\SW_X(\fs)$ and a universal
polynomial in the intersection form $Q_X$ and the classes
$c_1(\fs)-\Lambda$ and $\Lambda$, with coefficients depending only on
$(c_1(\fs)-\Lambda)^2$, $(c_1(\fs)-\Lambda)\cdot\Lambda$, $\Lambda^2$,
$\SW_X(\fs)$, and a universal function of $\chi(X)$ and $\sigma(X)$; the
polynomial degree depends on $\ell$ and the degree of the Donaldson
invariant. (In the paragraph following the statement of this conjecture in
\S \ref{subsec:MultiplicityConjecture} we give an informal explanation of
why Conjecture \ref{conj:Divisibility} should hold and in the last
paragraphs of \S \ref{subsec:MultiplicityConjecture} and \S
\ref{subsec:Compare}, we explain the role of the conjecture in the proof of
Theorem \ref{main}.) The work of the first and third authors (\cite{FL1},
\cite{FL2}, \cite{FL2a}, \cite{FL2b}, \cite{FLGeorgia},
\cite{FeehanGenericMetric}, \cite{FL3}, \cite{FL4}, \cite{FLKM1}) goes a
long way towards a proof of \cite[Conjecture 4.1]{FLGeorgia} and, in
particular, Conjecture \ref{conj:Divisibility}. If Conjecture
\ref{conj:Divisibility} is not assumed, a somewhat weaker result with a
correction term can still be deduced from those papers; the correction
depends on the maximum dimension of the non-empty Seiberg-Witten moduli
spaces --- see \S \ref{sec:FLthmDiscussed} for more details.  The abundance
condition is used to construct the right $\PU(2)$-monopole setup and is not
the most general condition under which the theorem can be proved; see the
discussion after Lemma~\ref{abunlem}.

If one assumes Witten's \cite{Witten} conjectured formula \eqref{wi}
relating the Seiberg-Witten invariants for manifolds of simple type and the
Donaldson invariants, the conclusion of Theorem~\ref{main} is equivalent to
a vanishing theorem for the Donaldson invariants through certain degrees.
A restricted version of Witten's conjecture --- see Theorem \ref{thm:FLthm}
below --- is proved in \cite{FL2}, \cite{FL2a}, \cite{FL2b}.  In fact this
theorem implies many different relations between the Donaldson and
Seiberg-Witten invariants and will be used here in two ways: to prove the
vanishing result for Donaldson invariants, Theorem~\ref{thm:DVanish} below,
and to deduce Theorem~\ref{main}, most cases of which follow from
Theorem~\ref{thm:DVanish} and the restricted version of Witten's
conjecture.

\begin{thm}
\label{thm:DVanish}
Assume Conjecture \ref{conj:Divisibility} and that $X$ is of SW-simple type
and abundant. Then for
$$
0\le d\le  c(X)-1\quad \mathrm{and}\quad m\ge 0,
$$
we have for all $w\in H^{2}(X,\ZZ)$ with $w= w_{2}(X) \pmod 2$,
$$
	D^{w}_{X}(h^{d-2m}x^{m})= 0.
$$
(See \S 2 for conventions concerning the Donaldson invariants $D^{w}_X$.)
\end{thm}

\subsubsection*{Acknowledgments}
We are grateful to Robert Friedman, Peter Ozsv\'ath, and Zolt\'an Sz\'abo for
helpful comments. The first author would like to thank the Institute for
Advanced Study, Princeton, the Institut des Hautes \'Etudes Scientifiques,
Bures-sur-Yvette, as well as the National Science Foundation, for their
generous support during the preparation of this article.

\section{The Donaldson invariants and relations to
the Seiberg-Witten invariants}

We quickly review some basic definitions regarding the Donaldson
invariants (see \cite{KMStructure}). Let
$$
\AAA(X) = \Sym\left( H_{\even}(X;\RR)\right)
$$
be the graded algebra, with $z=\beta_1\beta_2\cdots\beta_r$ having total
degree $\deg(z) = \sum_p(4-i_p)$, when $\beta_p\in H_{i_p}(X;\RR)$. In
particular a point $x \in X$ gives a distinguished generator still called
$x$ in $\AAA(X)$ of degree four. For any choice of $w\in H^{2}(X,\ZZ)$
there is a corresponding Donaldson invariant which is now a linear function
$$
D^{w}_{X}:\AAA(X) \to \RR
$$
and is defined by evaluating cohomology classes corresponding
to elements of $\AAA(X)$ on instanton moduli spaces
of $\SO(3)$-bundles $P$ with $w_{2}(P)\equiv w \pmod 2$, $w$
determining the orientation of the moduli spaces using Donaldson's conventions
\cite[\S 2(ii)]{KMStructure}.  If
$w\equiv w' \pmod 2$ we have
$$
D^{w}_{X}=(-1)^{\quarter (w-w')^{2}}D^{w'}_{X}.
$$
Also
$D^{w}_{X}(z)=0$ unless $z$ contains a monomial $m$ with
\begin{equation}
	\mathrm{deg}(m)\equiv -2w^{2}-\frac{3}{2}(\chi(X)+\sigma(X))\pmod 8,
	\label{mod8}
\end{equation}
and $\chi(X)+\sigma(X) \equiv 0 \pmod 4$.

\noindent
A four-manifold has KM-simple type if for all $z\in \AAA(X)$
we have
$$
D^{w}_{X}(x^{2}z)=4D^{w}_{X}(z).
$$
It is known that if this relation holds for one $w$, it holds for all
$w$.  For manifolds of KM-simple type one introduces the formal power
series in a variable $h \in H_{2}(X)$,
$$
\bD^{w}_{X}(h) = D^{w}_{X}((1+\half x)e^{h}).
$$
By Equation~\eqref{mod8} the series $\bD^{w}_{X}$ is an even function if
$$
-w^{2}+\frac{3}{4}(\chi(X)+\sigma(X))\equiv  0 \pmod 2,
$$
and is odd otherwise. Notice that $\bS\bW_X^{w}$ has the same
property since
$$
SW_X(\fs)=
(-1)^{\quarter(\chi(X)+\sigma(X))}\SW_X(\bar\fs).
$$
In addition if $w\equiv w' \pmod 2$ then
$$
\bS\bW_X^{w'}=(-1)^{\quarter (w-w')^{2}}\bS\bW_X^{w}.
$$

\noindent
According to \cite[Theorem 1.7]{KMStructure}, when $X$ has KM-simple type
the series $\bD^{w}_{X}(h)$ is an analytic function of $h$ and there are
finitely many characteristic cohomology classes $K_{1},\ldots,K_{m}$ (the
KM-basic classes) and constants $a_{1},\ldots,a_{m}$ (independent of $w$)
so that
$$
\bD^{w}_X(h)= e^{\half h \cdot h}\sum_{i=1}^{r}(-1)^{\half(w^{2}+K_{i}\cdot
w)}a_{i}e^{\langle
K_{i},h\rangle}.
$$
Witten's conjecture \cite{Witten} relating the Donaldson and Seiberg-Witten
invariants of manifolds of simple type says that
\begin{equation}\label{wi}
	\bD^{w}_{X}(h)=2^{2-c(X)}e^{\half h\cdot h}\bS\bW^{w}_{X}(h).
\end{equation}

While a complete, mathematically rigorous proof of Witten's conjecture has
not yet been obtained, a possible approach using a $\PU(2)$-monopole
cobordism was proposed by Pidstrigach and Tyurin \cite{PTLocal}.  By
employing the $\PU(2)$-monopole cobordism, the first and third author
\cite{FL2}, \cite{FL2a}, \cite{FL2b} proved relations between these
invariants which we restate here in somewhat restricted form.  To state
their results define for $\Lambda\in H^{2}(X;\ZZ)$,
\begin{equation}
\label{eq:SWSimpleUDepthParam}
r(\Lambda) =  -\Lambda^{2}-\frac{1}{4}(11\chi(X)+15\sigma(X)).
\end{equation}
This number keeps track of the depth in  the Uhlenbeck
compactification of the $\PU(2)$-moduli space that reducibles appear
in terms of the degree of the corresponding Donaldson invariant when $X$
has SW-simple type.
Also define
\begin{equation}
\label{eq:DiracIndexParam}
i(\Lambda) = \Lambda^{2}-\frac{1}{4}(3\chi(X)+7\sigma(X)).
\end{equation}
This number keeps track of the index of the Dirac operator for the
$\PU(2)$-moduli space in terms of the degree of the corresponding Donaldson
invariant. If we assume Conjecture~\ref{conj:Divisibility} the results of
\cite{FL2}, \cite{FL2a}, \cite{FL2b} yield:

\begin{thm}\label{thm:FLthm}
\cite[Theorem 1.1]{FL2b} Assume that $X$ satisfies
	Conjecture~\ref{conj:Divisibility} and is of SW-simple type.
	Suppose that $\Lambda \in B^{\bot}$ and that $w \in H^{2}(X;\ZZ)$ with
	$w-\Lambda\equiv w_{2}(X) \pmod 2$.  If we have
	$$
	\delta< r(\Lambda)\quad \mathrm{and} \quad  \delta <i(\Lambda),
	$$
	then for all $h\in H^2(X;\RR)$,
	\begin{equation}\label{eq:Dvanish}
		D^{w}_{X}(h^{\delta-2m}x^{m})= 0.
	\end{equation}
	
	\noindent
	If we instead have
	$$
	\delta= r(\Lambda)\quad \mathrm{and} \quad  \delta <i(\Lambda),
	$$
	then for all $h\in H^2(X;\RR)$,
	\begin{equation}
	\label{eq:DSWrel}
\begin{aligned}
	D^{w}_{X}(h^{\delta-2m}x^{m})
&=
	2^{1-\half (c(X)+\delta)}(-1)^{m-1+\half\Lambda^2-\Lambda\cdot w}
\\
&\qquad\times \sum_{\fs \in  S}
	(-1)^{\half(w^{2}+c_{1}(\fs)\cdot w)}SW_X(\fs)
	\langle c_{1}(\fs)-\Lambda,h \rangle^{\delta -2m}.
\end{aligned}
	\end{equation}
	\end{thm}

Since $\Lambda^2-2\Lambda\cdot w = \sigma(X) - w^2 \pmod{8}$, the sign
factor $(-1)^{\half\Lambda^2-\Lambda\cdot w}$ may also be written as
$(-1)^{\half(\sigma(X) - w^2)}$.

The hypotheses of Theorem 1.1 in \cite{FL2b} do not require that
Conjecture~\ref{conj:Divisibility} holds, that $b_2^+(X)>1$, or that $X$
have Seiberg-Witten simple type.  If $b_2^{+}(X)=1$ a similar result also
holds when the chamber structure is taken into account. See \cite[Theorem
1.1]{FL2b} for the most general statements, including a partial treatment
of the case when $b_{1}(X)>0$.  The role of SW-simple type, the manner in
which Conjecture~\ref{conj:Divisibility} strengthens Theorem 1.1 in
\cite{FL2b} and the derivation of Theorem \ref{thm:FLthm} from it are
explained in \S \ref{sec:FLthmDiscussed}. The mechanism underlying the
vanishing result in Equation~\eqref{eq:Dvanish} was pointed out in
\cite{FL2}: see, for example, Theorem 5.33 and the paragraph following
Equation~(7.2). The fact that the $\PU(2)$-monopole cobordism should imply
some relations between the SW-basic classes was evident from Lemmas 5.30
and 5.31 in \cite{FL2}: see Conjecture 4.1 in \cite{FLGeorgia}.

	Elements of $B$ are always characteristic, so the
	condition $\Lambda\in B^{\bot}$ (assuming $B$ is non-empty) means
	that $\Lambda\cdot w_{2}$ is zero, and $\Lambda^{2}$ is therefore
	even.  Via the additional condition $\Lambda - w \equiv w_{2}\pmod
	2$, the pair $w$ and $\Lambda$ is constrained so that $w^{2} \equiv
	\Lambda^{2}+\sigma(X) \pmod 4$.  Furthermore, if $w$ is
	characteristic, so $w^2\equiv \sigma(X)\pmod 8$, and $B$ is
	non-empty then this condition implies that $\Lambda\equiv 0\pmod 2$
	and $\Lambda^{2}\equiv 0 \pmod 8$.

	Figure 1 below shows a typical picture in the
	$(\Lambda^{2},\delta)$ plane, assuming that $c(X)$ is positive.
		The
	lines $\delta=i(\Lambda)$ and $\delta=r(\Lambda)$ intersect at
	$\delta=c(X)$ and $\Lambda^{2}= - (\chi(X)+\sigma(X))$.
		The Theorem asserts the vanishing of a Donaldson invariant
when
	$(\Lambda^{2},\delta)$ is in the interior of the outlined triangle,
	and gives a formula for the Donaldson invariant in terms of
	Seiberg-Witten invariants when $(\Lambda^{2},\delta)$ is one of the
	marked points on the right-hand edge of the triangle.

	The marked points are the lattice points with $2\delta\equiv
	-2w^{2}-\frac{3}{2}(\chi(X)+\sigma(X))\pmod 8 $ (see
	Equation~\eqref{mod8}) and $\Lambda^{2}\equiv w^{2}-\sigma(X) \pmod
	4$. We have drawn the case in which $w^{2}-\sigma(X)$ is zero mod
	$4$, so that the intersection of the two lines is one of the marked
	lattice points (because $2c(X)$ always satisfies the condition
	\eqref{mod8}).  In the case that $w$ is characteristic,
	$\Lambda^{2}$ is constrained to be $0$ mod $8$ (the white dots); we
	have drawn the case that $-(\chi(X)+\sigma(X))$ is $0$ mod $8$
	also.

\begin{figure}
    \begin{center}
		\begin{picture}(200,270)(80,0)
		\multiput (135 ,215)(20,-20){9}{\line(1,-1){10}}
		\multiput (55,5)(20,20){11}{\line(1,1){10}}
		\put (35,80){\vector(1,0){275}}
		\put (80,0){\vector(0,1){225}}
		\put (315,75){$\Lambda^{2}$}
		\put (85,225){$\delta$}
		\put (130,225){$\delta=r(\Lambda)$}
		\put (270,225){$\delta=i(\Lambda)$}
		\put (25,152){$\delta=c(X)$}
		\put (75,150) {\line(1,0){10}}
		\put (150,55) {$\Lambda^{2}=-(\chi(X)+\sigma(X))$}
		\put (200,75) {\line(0,1){10}}
		\multiput(40,10)(0,20){11}{\circle{2}}
		\multiput(60,10)(0,20){11}{\circle*{2}}
		\multiput(80,10)(0,20){11}{\circle{2}}
		\multiput(100,10)(0,20){11}{\circle*{2}}
		\multiput(120,10)(0,20){11}{\circle{2}}
		\multiput(140,10)(0,20){11}{\circle*{2}}
		\multiput(160,10)(0,20){11}{\circle{2}}
		\multiput(220,10)(0,20){11}{\circle*{2}}
		\multiput(240,10)(0,20){11}{\circle{2}}
		\multiput(260,10)(0,20){11}{\circle*{2}}
		\multiput(280,10)(0,20){11}{\circle{2}}
		\multiput(300,10)(0,20){11}{\circle*{2}}
		\multiput(200,10)(0,20){11}{\circle{2}}
		\multiput(180,10)(0,20){11}{\circle*{2}}
		\qbezier[10](160,85)(150,85)(160,95)
		\qbezier[20](160,95)(180,115)(195,130)
		\qbezier[10](195,130)(200,135)(205,130)
		\qbezier[20](205,130)(220,115)(240,95)
		\qbezier[10](240,95)(250,85)(240,85)
		\qbezier[30](240,85)(200,85)(160,85)
		\qbezier[20](215,130)(230,115)(255,90)
		\qbezier[20](225,130)(240,115)(265,90)
		\qbezier[10](255,90)(265,85)(265,90)
		\qbezier[10](215,130)(215,135)(225,130)
		\put (67,70) {$O$}
	\end{picture}
  \end{center}
	  \caption{}
	  \label{triangle}
\end{figure}

	Thus we see that in order to apply Theorem~\ref{thm:FLthm} we need
	the restriction of the intersection form to $B^\perp$ to be rich
	enough that it realizes particular values, on classes with certain
	constraints on their mod 2 reductions.  The definition of
	abundant is chosen simply to have some rather general but compact
	hypothesis under which we can find the required classes.  The
	precise hypothesis needed are the conclusions of the following
	Lemma.

\begin{lem}\label{abunlem} If $X$ is abundant then there are
	cohomology classes $\Lambda_{0},\Lambda_{1}\in B^{\bot}$ with
	 $ \Lambda_{0}\equiv \Lambda_{1}\pmod 2$
	so that:
	       $$
            \Lambda^{2}_{0}=-(\chi(X)+\sigma(X)),
	         \quad \Lambda_{1}^{2}=-(\chi(X)+\sigma(X))+4.
	      $$
        There is also a class $\Lambda \in 2 B^{\bot}$ with
	$\Lambda^{2}=-(\chi(X)+\sigma(X))$ if $-(\chi(X)+\sigma(X))\equiv
	0\pmod 8$ and $\Lambda^{2}=-(\chi(X)+\sigma(X))+4$ if
	$-(\chi(X)+\sigma(X))\equiv
	4\pmod 8$.
\end{lem}
\begin{proof}
     Since $X$ is abundant we can find $e_{1}, e_{2} \in B^{\bot}$
	 so that $e_{1}\cdot e_{1} = e_{2}\cdot e_{2}=0$ and $e_{1}\cdot
	 e_{2}=1.$   Set $h=\quarter(\chi(X)+\sigma(X))$.
	 Then $\Lambda_{0}= 2e_{1}-he_{2}$ and
	 $\Lambda_{1}=2e_{1}+(1-h)e_{2}$ will do. If $h\equiv 0\pmod 2$
	 then taking $\Lambda= \Lambda_{0}$ proves the second assertion
	 and if $h\equiv 1 \pmod 2$ taking $\Lambda=\Lambda_{1}$
	 proves the second assertion.
\end{proof}
There are non-abundant four-manifolds: for example, some of
the fake K3-surfaces of \cite{GompfMrowka} fail to be abundant.
If log transforms are performed on tori in three
distinct nuclei then the intersection form on $B^{\bot}$ is a degenerate
form with three-dimensional radical and having an $-E_{8}\oplus
-E_{8}$ summand.  These manifolds
however satisfy the conclusion of Lemma~\ref{abunlem}.
On the other hand one can show that a simply connected minimal
surface of general type is abundant.

Here is the proof of Theorem~\ref{main} of the introduction.

\begin{proof}[Proof of Theorem~\ref{main}]
	Choose $w$  characteristic as in the hypotheses of Theorem~\ref{main}.
	To prove the theorem we must show that the Taylor
	coefficients of $\bS\bW_X^{w}$ are zero
	in degrees $d$ with $0\le d \le c(X)-3$.
	These coefficients are zero unless $d\equiv
	-w^{2}+\frac{3}{4}(\chi(X)+\sigma(X)) \equiv c(X) \pmod
	2$ so write
	$d= c(X)-4-2m$ for some $m\ge 0$.  Take
	$\Lambda_{0}, \Lambda_{1} \in B^{\bot}$ with
	$$
	\Lambda_{0}^{2}= -(\chi(X))+\sigma(X)),\quad
	\Lambda_{1}^{2}=-(\chi(X)+\sigma(X))+4
	$$
	and $\Lambda_{0}\equiv \Lambda_{1} \pmod 2$ as provided
	by Lemma~\ref{abunlem}.
	We have
	$$
	       r(\Lambda_{0})=c(X)=
		   i(\Lambda_{0})
	$$
	so taking $\delta =d+2m=c(X)-4$
	in Equation~\eqref{eq:Dvanish}
	of Theorem~\ref{thm:FLthm} 
(applied to calculate $D^{w'}_X$, with $w'=w+\Lambda_0$ and
	$\Lambda=\Lambda_0$), one concludes that
          $$
          D^{w+\Lambda_{0}}_{X}(h^{d}x^{m})=0.
	$$
	On the other hand we have
	$$
	r(\Lambda_{1})=c(X)-4 \quad \mathrm{and}\quad
	i(\Lambda_{1})=c(X)+4
	$$
	so taking $\delta= c(X)-4$ again but now using
Equation~\eqref{eq:DSWrel} (applied to calculate $D^{w'}_X$, with
$w'=w+\Lambda_1$ and $\Lambda=\Lambda_1$), one concludes that
	\begin{equation}
			\begin{aligned}
			D^{w+\Lambda_{1}}_{X}(h^{d}x^{m})
			&=
			2^{1-\half
(c(X)+d)-m}(-1)^{m-1+\half\Lambda_1^2-\Lambda_{1}\cdot w}
			\\
			&\qquad\times \sum_{\fs \in  S}
			(-1)^{\half(w^{2}+c_{1}(\fs)\cdot w)}SW_X(\fs)
			\langle c_{1}(\fs)-\Lambda_{1},h \rangle^{d}.
		         \end{aligned}
			\label{eq:1}
	\end{equation}
	Since $\Lambda_{0}\equiv \Lambda_{1}\pmod 2$ and
	$\Lambda_{0}-\Lambda_{1}\in B^{\bot}$
	we have
	$0=D_{X}^{w+\Lambda_{0}}(h^{d}x^{m})=
	D_{X}^{w+\Lambda_{1}}(h^{d}x^{m})$.
	Putting all this together implies that the
	Taylor coefficients at the origin of the analytic function
          $$
             e^{-\langle\Lambda_{1},h\rangle}\bS\bW_{X}^{w}(h)
          $$
         are zero up to degree $c(X)-3$. Since
          $e^{-\langle\Lambda_{1},h\rangle}$ is invertible
	it must be that the $\bS\bW_{X}^{w}(h)$ vanishes to
         the required order there.
\end{proof}

\noindent
It remains to prove Theorem~\ref{thm:DVanish} of the introduction.
\begin{proof}[Proof of Theorem~\ref{thm:DVanish}]
	Fix $w$ characteristic.  Now there is nothing to prove unless
	$d\equiv
	-w^{2}+\frac{3}{4}(\chi(X)+\sigma(X)) \pmod
	4$ so assume this as well.  If $-(\chi(X)+\sigma(X))\equiv 0 \pmod
	8$ find $\Lambda \in 2B^{\bot}$ so that
	$\Lambda^{2}=-(\chi(X)+\sigma(X))$ as in Lemma~\ref{abunlem}.  Then
	$$
	i(\Lambda)=r(\Lambda)=c(X)
	$$
	and so for any $d\le c(X)-4$ we can use Equation~\eqref{eq:Dvanish}
	to conclude
	$$
	0=D_{X}^{w+\Lambda}(h^{d-2m}x^{m})= D_{X}^{w}(h^{d-2m}x^{m}).
	$$
	If $-(\chi(X)+\sigma(X))\equiv 4\pmod 8$ find $\Lambda \in
	2B^{\bot}$ so that
	$\Lambda^{2}=-(\chi(X)+\sigma(X))+4$ as in Lemma~\ref{abunlem}.  Then
	$$
	i(\Lambda)=c(X)+4\quad\mathrm{and}\quad r(\Lambda)=c(X)-4
	$$
	and so for any $d\le c(X)-8$ we can use Equation~\eqref{eq:Dvanish}
	to conclude
	$$
	0 = D_{X}^{w+\Lambda}(h^{d-2m}x^{m})= D_{X}^{w}(h^{d-2m}x^{m}).
	$$
	While for $d=c(X)-4$ Equation~\eqref{eq:DSWrel} applies to give
	\begin{equation}
			\begin{aligned}
			D^{w}_{X}(h^{d-2m}x^{m})
			&=
			2^{1-\half (c(X)+d)}
(-1)^{m+\half\Lambda^2-\Lambda\cdot w}
			\\
			&\qquad\times \sum_{\fs \in  S}
			(-1)^{\half(w^{2}+c_{1}(\fs)\cdot w)}SW_X(\fs)
			\langle c_{1}(\fs)-\Lambda,h \rangle^{d -2m}.
		         \end{aligned}
			\label{eq:2}
	\end{equation}
	But up to a constant factor the right-hand side of this equation is
	the Taylor coefficient of degree $d-2m$ in $h$ of
	$e^{-\langle\Lambda,h\rangle}\bS\bW_{X}^{w}(h)$
	and hence vanishes by the previous theorem.
\end{proof}

Notice that it is only the case $-(\chi(X)+\sigma(X))\equiv 4\pmod 8$ and
$d=c(X)-4$ which requires the formula \eqref{eq:DSWrel} (and thus
\cite[Theorem 1.1]{FL2b}); the remaining cases are direct consequences of
the more elementary vanishing result \eqref{eq:Dvanish} (and thus
\cite[Theorem 3.33]{FL2b}).  Also notice that
in many cases the vanishing can be deduced from another mechanism.  If
there are distinct classes $\Lambda$ and $\Lambda'$ so that
Equation~\eqref{eq:2} holds then we deduce that
$e^{-\langle\Lambda,h\rangle}\bS\bW_{X}^{w}(h)
=e^{-\langle\Lambda',h\rangle}\bS\bW_{X}^{w}(h)$ up to order $c(X)-4$ in
$h$ and hence $\bS\bW_{X}^{w}(h)$ vanishes to that order.


\section{On Theorem \ref{thm:FLthm} and Conjecture \ref{conj:Divisibility}}
\label{sec:FLthmDiscussed}
We explain how to derive Theorem \ref{thm:FLthm} from \cite[Theorem
1.1]{FL2b} and Conjecture \ref{conj:Divisibility}.

\subsection{The multiplicity conjecture}
\label{subsec:MultiplicityConjecture}
The formula in \cite[Theorem 1.1]{FL2b} for Donaldson invariants in terms
of Seiberg-Witten invariants requires a choice of base \spinc structure
$\fs_0=(W^+,W^-,\rho)$ and Hermitian, rank-two vector bundle $E$ so that
reducible $\PU(2)$ monopoles (which we view as Seiberg-Witten monopoles ---
see \cite[\S 2 \& 3]{FL2a}) appear only in the top level of $\barM_{W,E}$,
if at all.  However, by assuming Conjecture \ref{conj:Divisibility}, we can
significantly relax the preceding constraint and allow reducible $\PU(2)$
monopoles to appear in the lower-levels of $\barM_{W,E}$ provided their
associated Seiberg-Witten invariants vanish. Hence, the sum in
\cite[Theorem 1.1]{FL2b} over set of \spinc structures defining reducible
$\PU(2)$ monopoles lying only in the top level $M_{W,E}$ can be replaced by
a sum over \spinc structures whose associated reducibles can lie in the top
level of $\barM_{W,E}$ if their Seiberg-Witten invariants are non-zero and
in any level if they are zero. This is a much weaker requirement than that
of \cite[Theorem 1.1]{FL2b}, which asks that there be {\em no non-empty\/}
lower-level Seiberg-Witten moduli spaces, a condition which depends, for
example, on the choice of perturbations and is difficult to verify in
practice.

By the hypothesis of Theorem \ref{thm:FLthm}, we are given a class $w\in
H^2(X;\ZZ)$: recall from Equation~\eqref{mod8} that the Donaldson invariant
$D_X^w(z)$ of $\deg(z) = 2\delta\in 2\ZZ$ is defined to be zero unless
$\delta \equiv -w^2 -\frac{3}{4}(\chi+\sigma) \pmod{4}$.  As there is no
loss in assuming this, let $p\in\ZZ$ be determined by the equation
\begin{equation}
\label{eq:Choiceofdelta}
\delta = -p - \frac{3}{4}(\chi+\sigma),
\end{equation}
so $p = w^2 \pmod{4}$, and define a Hermitian, rank-two vector bundle $E$ over
$X$ by requiring that $c_1(E) = w$ and $c_2(E)-\quarter c_1(E)^2 =
-\quarter p$. For such a bundle $E$, we have $p_1(\su(E)) = p$. 

The hypotheses of Theorem \ref{thm:FLthm} also provide us with a class
$\Lambda \in H^2(X;\ZZ)$. We therefore fix a \spinc structure $\fs_0 =
(W^+,W^-,\rho)$ on $X$ with $c_1(\fs_0) = c_1(W^+)$ determined by
$$
c_1(W^+) + c_1(E) = \Lambda.
$$
Any other \spinc structure $\fs$ on $X$ defines a class in $H^2(X;\ZZ)$ and
a Hermitian line bundle $L_1$ over $X$ such that $c_1(\fs) = c_1(\fs_0) +
2c_1(L_1)$. We write $\fs = \fs_0\otimes L_1 = (W^+\otimes L_1,W^-\otimes
L_1,\rho)$. Hence, the \spinc structure $\fs$ determines a split, rank-two,
Hermitian vector bundle $E' = L_1\oplus (\det E)\otimes L_1^*$ with
$p_1(\su(E')) = p_1(\su(E)) \pmod{4}$ and
\begin{equation}
\label{eq:Pontrjagin}
p_1(\su(E')) = (2c_1(L_1) - c_1(E))^2= (c_1(\fs)- \La)^2. 
\end{equation}
These split bundles, $E'$, in turn define families of reducible $\PU(2)$
monopoles.

The vector bundle $E'$ defines reducible $\PU(2)$ monopoles in some level of
$\barM_{W,E}$ only if $p_1(\su(E'))\geq p_1(\su(E))$, in which case they
would be contained in the level $\ell\in\ZZ_{\geq 0}$, where 
\begin{equation}
\label{eq:Level}
p_1(\su(E'))
= 
p_1(\su(E))) + 4\ell.
\end{equation}
Let us first consider the case $\ell>0$, so the corresponding reducible
solutions to the $\PU(2)$ monopole equations lie in the level
$M_{W,E_{-\ell}}\times \Sym^\ell(X)$ of the Uhlenbeck compactification
$$
\barM_{W,E}
\subset
\bigcup_{\ell=0}^N M_{W,E_{-\ell}}\times \Sym^\ell(X),
$$
where $E_{-\ell}$ is a Hermitian, rank-two vector bundle over $X$ with
$c_1(E_\ell)=c_1(E)$ and $c_2(E_{-\ell})=c_2(E)-\ell$.  The smooth loci
$M^{*,0}_{W,E_{-\ell}}\times \Si$ of the lower level $M_{W,E_{-\ell}}\times
\Sym^\ell(X)$, where $\Si\subset\Sym^\ell(X)$ is a smooth stratum, have
codimension greater than or equal to $2\ell$.  For $z\in\AAA(X)$, a
geometric representative $\barV(z)\subset \barM_{W,E}$ of codimension
$\deg(z)$ is defined in \cite{FL2b}.  It is shown in \cite{FL2b} that this
geometric representative intersects the lower levels of $\barM_{W,E}$ in a
set with codimension $\deg(z)$ except at the reducible points.  There is
also a geometric representative $\barW(x^{n_{c_1}})\subset \barM_{W,E}$
which has codimension $2n_{c_1}$ on the complement of the reducible and
zero-section points.  The formula in \cite[Theorem 3.33]{FL2b} is proved by
computing the intersections of the geometric representatives,
$\barV(z)\cap\barW(x^{n_{c_1}})$, with the links of the strata of the
anti-self-dual and reducible $\PU(2)$ monopoles in $M_{W,E}$.  Thus to
extend \cite[Theorem 3.33]{FL2b} to the case where there are reducibles in
the lower levels $M_{W,E_{-\ell}}\times \Sym^\ell(X)$, $\ell>0$, one needs
to consider the links in $\barM_{W,E}$ of the families
\begin{equation}
\label{eq:LowerLevelReducible}
M^{sw}_\fs\times\Sym^\ell(X).
\end{equation}
Here, $M^{sw}_{\fs}$ is the perturbation of the
standard Seiberg-Witten moduli space described in \cite[Equation
(2.14)]{FL2a} and \cite[Lemma 3.13]{FL2a}.

The gluing theorems of \cite{FL3}, \cite{FL4}
provide a sufficiently explicit
description of a neighborhood of the family \eqref{eq:LowerLevelReducible}
that we can define their link in $\barM_{W,E}$, which we denote by
$\bL_{W,E,L_1}$.  Given the results of \cite{FL3}, \cite{FL4}, the formula of
\cite[Theorem 3.33]{FL2b} can then be replaced by
\begin{equation}
\label{eq:ConjectureSum}
D^w_X(z) = -2^{1-n_a}\sum_{\fs\in\Spinc(X)}
\# 
\left( 
\barV(z)\cap\barW(x^{n_a-1}) \cap \bL_{W,E,L_1}
\right).
\end{equation}
This is a sum over the finite set of \spinc structures $\fs = \fs_0\otimes
L_1$ with a non-empty SW moduli space $M^{sw}_\fs$, defining reducible
$\PU(2)$ monopoles contained in $\barM_{W,E}$.  The definition of the links
$\bL_{W,E,L_1}$ and the proof that the intersection numbers in
Equation~\eqref{eq:ConjectureSum} are well-defined appears in
\cite{FL4}. If we assume the following conjecture, the sum in
Equation~\eqref{eq:ConjectureSum} can be reduced to a sum over \spinc
structures $\fs$ with $\SW_X(\fs)\neq 0$.

\begin{conj}
\label{conj:Divisibility}
Let $X$ be a closed, oriented, smooth four-manifold with $b_2^+(X)>0$.  Let
$\fs=\fs_0\otimes L_1$ be a \spinc structure defined by $\fs_0$ and a
reduction $E_{-\ell}=L_1\oplus (\det E)\otimes L_1^*$, where
$c_2(E_{-\ell})=c_2(E)-\ell$ and $\ell\geq 0$.  Let $\bL_{W,E,L_1}$ be the
link of the ideal Seiberg-Witten moduli space $M^{sw}_\fs\times
\Sym^\ell(X)$. Then the intersection number
$$
\#
\left(
\barV(z)\cap\barW(x^{n_a-1}) \cap \bL_{W,E,L_1}
\right)
$$
appearing in the right-hand side of Equation~\eqref{eq:ConjectureSum} is a
multiple of $\SW_X(\fs)$ and thus vanishes if $\SW_X(\fs)=0$.
\end{conj}

For reducible $\PU(2)$ monopoles contained in the top level, Conjecture
\ref{conj:Divisibility} follows immediately from Theorem 1.1 in \cite{FL2b}
and, when $\ell=1$, from unpublished work of the first and third author
\cite{FLinprep}. To appreciate why the conjecture should hold more
generally, recall that the construction in \cite{FL3} of the link
$\bL_{W,E,L_1}$ via gluing maps shows that $\bL_{W,E,L_1}$ is a union of
pieces $\bL_\Si$, where $\Sigma$ ranges over the smooth strata of
$\Sym^\ell(X)$. Each $\bL_\Si$ admits a fiber-bundle structure over
$M^{sw}_{\fs}\times \Si$. By choosing suitable lifts of the cohomology
classes on the link $\bL_{W,E,L_1}$ to cohomology classes with compact
support in the pieces $\bL_\Si$, a pushforward argument shows that each
pairing with $\bL_\Si$ can be expressed in terms of a pairing with
$M^{sw}_{\fs}\times \Sigma$.  The cohomology classes (given by the
pushforward) to be paired with $M^{sw}_{\fs}\times\Si$ all contain a factor
of $h^{d_s}$, where $2d_s=\dim M^{sw}_{\fs}$ and $h$ is the element of
$H^2(M^{sw}_{\fs};\ZZ)$ used to define the Seiberg-Witten invariant.  Thus, all
the pairings would be multiples of $\SW_X(\fs)$.  This argument is
discussed further in \cite[pp. 141--142]{FLGeorgia} for the case of a
zero-dimensional moduli space $M^{sw}_{\fs}$.

If we did not assume the conjecture, the conditions on $\Lambda$ necessary to
ensure that the reducible points in $\barM_{W,E}$ appear only in the top
level would be much stricter than those of Theorem \ref{thm:FLthm} (and
generally unverifiable).  Not only would we need to require that
$\La\cdot c_1(\fs)=0$ for all $\fs$ with $M^{sw}_{\fs}$ non-empty (rather
than just $\fs$ with $\SW_X(\fs)\neq 0$), but the presence of
positive-dimensional, non-empty Seiberg-Witten moduli spaces would further
reduce the degree of Donaldson invariant computable by Theorem
\ref{thm:FLthm} (see the explanation of the condition $\delta\le r(\La)$ in
\S \ref{subsec:Compare} below).  Finally, moduli spaces $M^{sw}_{\fs}$
with $\SW_X(\fs)=0$ may or may not be empty, depending on the choice of
perturbations, and --- except for some special cases \cite{Witten} --- it
appears very difficult to choose perturbations in such a way that those
moduli spaces would be empty; moreover, the perturbations of the
Seiberg-Witten equations arising as reductions of $\PU(2)$ monopole
equations are constrained by the choice of perturbations in the latter
equations \cite{FeehanGenericMetric}, \cite{FL2a}, \cite{FL2b}.

\subsection{Comparing indices and degrees}
\label{subsec:Compare}
Continue the notation of \S \ref{subsec:MultiplicityConjecture}.
To apply Theorem 1.1 in \cite{FL2b}, given Conjecture
\ref{conj:Divisibility}, we need to check that the following conditions are
obeyed:
\begin{itemize}
\item 
The \spinc structures $\fs$ with non-zero Seiberg-Witten invariants
define reducible $\PU(2)$ monopoles which are contained only in
the top level of $\barM_{W,E}$, if at all, and
\item 
The index $n_a = \Ind_\CC D_A$ is positive, where $D_A:\Gamma(W^+\otimes
E) \to \Gamma(W^-\otimes E)$ is the Dirac operator appearing in the 
$\PU(2)$ monopole equations.
\end{itemize}
We now show how these conditions are equivalent to the constraints on
$\delta$, $r(\La)$ and $i(\La)$ in the hypotheses of Theorem
\ref{thm:FLthm}.

\begin{proof}[Proof that Theorem 1.1 in \cite{FL2b} 
and Conjecture \ref{conj:Divisibility} implies Theorem \ref{thm:FLthm}] For
$\fs\in S$, let $\fs$ (and the choice of $W^+$) define a Hermitian,
rank-two vector bundle $E'$ with $c_1(E')=c_1(E)$, as in \S
\ref{subsec:MultiplicityConjecture}.  We have $c_1(\fs)^2 = 2\chi(X)
+3\si(X)$ by Equation~\eqref{eq:SimpleType} (since $X$ has SW-simple type)
and $c_1(\fs)\cdot \Lambda=0$ (by choice of $\Lambda$), so the formula
\eqref{eq:Pontrjagin} for $p_1(\su(E'))$ gives
$$
p_1(\su(E'))
=
(c_1(\fs)- \Lambda)^2
=
2\chi(X) +3\si(X) + \Lambda^2.
$$
By the definition \eqref{eq:Level} of $\ell$, the identity
\eqref{eq:Choiceofdelta} for $\delta$ (with $p = p_1(\su(E))$), and the
definition \eqref{eq:SWSimpleUDepthParam} of $r(\Lambda)$, we have
\begin{align*}
4\ell 
&= p_1(\su(E')) - p_1(\su(E))
\\
&= \delta  + \Lambda^2 + \frac{1}{4}\left( 11\chi(X) + 15\si(X)\right)
\\
&= \delta - r(\La).
\end{align*}
Therefore, $M^{sw}_{\fs}$ is contained in $\barM_{W,E}$ only if 
$\delta \geq r(\Lambda)$, in which case it lies in the level 
$$
\ell(\delta,\Lambda) 
=
\frac{1}{4}(\delta - r(\Lambda)).
$$
Note that in arriving at the formula $\ell(\delta,\Lambda)$, we only used
the facts that $c_1(\fs)^2=2\chi(X)+3\sigma(X)$ and
$c_1(\fs)\cdot\Lambda=0$, for $\fs\in S(X)$: thus if any one SW-basic class
$c_1(\fs)$ defines a reducible in the level $\ell$ of $\barM_{W,E}$, then
all the reducibles associated to SW-basic classes are contained in this
level.

The second condition necessary to apply \cite[Theorem 1.1]{FL2b} is that the
index of the Dirac operator, $D_A$, be positive.  Using the identity
\begin{equation}
\label{eq:IndexDirac}
\Ind_\CC D_A 
=
\frac{1}{4} (p_1(\su(E)) +\La^2 -\si(X)),
\end{equation}
and the identity \eqref{eq:Choiceofdelta} for $\delta$, together with the
definition \eqref{eq:DiracIndexParam} of the parameter
$i(\Lambda)$, one sees that
$$
\Ind_\CC D_A = \frac{1}{4}(i(\La)-\delta).
$$
Hence, the condition $\delta < i(\La)$ in the hypotheses of Theorem
\ref{thm:FLthm} is equivalent to $\Ind_\CC D_A>0$ 
and so we can apply Theorem 1.1 in \cite{FL2b}. 

Since $\delta \leq r(\Lambda)$, the \spinc structures $\fs$ defining basic
classes give Seiberg-Witten moduli spaces $M_\fs^{sw}$ and reducible
$\PU(2)$ monopoles contained only in the top level of $\barM_{W,E}$. The
only other possible contributions to the Donaldson invariants
$D_X^w(h^{\delta-2m}x^m)$ would be due to pairings of cohomology classes
with links of Seiberg-Witten moduli spaces $M_{\fs'}^{sw}$ defining
reducible $\PU(2)$ monopoles contained in levels $\ell\geq 1$ of
$\barM_{W,E}$. As such classes $c_1(\fs')$ are not basic, we have
$\SW_X(\fs')=0$ and so Conjecture \ref{conj:Divisibility} implies that the
moduli spaces $M_{\fs'}^{sw}$ make no contribution to the Donaldson
invariant $D_X^w(h^{\delta-2m}x^m)$.
\end{proof}

In Equation~\eqref{eq:Dvanish} of Theorem \ref{thm:FLthm}, which asserts
the vanishing of Donaldson invariants through certain degrees, one only
needs the more elementary Theorem 3.33 from \cite{FL2b} --- as it is not
necessary to evaluate the pairings on the right-hand-side --- rather than the
more difficult Theorem 1.1 from \cite{FL2b} which is needed to compute the
formula \eqref{eq:DSWrel} for Donaldson invariants.

Without the assumption that $X$ have SW-simple type, we could at best write
$$
c_1(\fs)^2 \leq 2\chi(X) + 3\sigma(X) + 4d(X),
$$
where $d(X)$ is the maximal expected dimension of the Seiberg-Witten moduli
spaces corresponding to $\fs \in S$. The results
of \S \ref{sec:Intro} could then be replaced with weaker ones, where the
order of vanishing would now also depend on $d(X)$.

In the absence of Conjecture \ref{conj:Divisibility}, we would need to
replace the set $S$ in the statement of Theorem
\ref{thm:FLthm} by the subset of $H^2(X;\ZZ)$ corresponding to \spinc
structures with non-empty Seiberg-Witten moduli spaces. Furthermore, the
SW-simple type equation for $c_1(\fs)^2$ would have to be replaced by the
weaker inequality described in the preceding remark, but with $d(X)$ now
denoting the maximal dimension of the non-empty Seiberg-Witten moduli
spaces. (In this case, $d(X)$ is no longer independent of the choice of
generic perturbation parameters.)


\ifx\undefined\bysame
\newcommand{\bysame}{\leavevmode\hbox to3em{\hrulefill}\,}
\fi

\end{document}